# GENERALIZED BOOTSTRAP FOR ESTIMATING EQUATIONS

By Snigdhansu Chatterjee and Arup Bose

*University of Minnesota and Indian Statistical Institute, Kolkata*

We introduce a generalized bootstrap technique for estimators obtained by solving estimating equations. Some special cases of this generalized bootstrap are the classical bootstrap of Efron, the deleted-$d$ jackknife and variations of the Bayesian bootstrap. The use of the proposed technique is discussed in some examples. Distributional consistency of the method is established and an asymptotic representation of the resampling variance estimator is obtained.

**1. Introduction.** One of the most popular ways of obtaining estimators for parameters in statistics is by solving "estimating equations." Examples are abundant in the contexts of quasi-likelihood methods, time series, biostatistics, stochastic processes, spatial statistics, robust inference, survey sampling and other areas. Godambe (1991) and Basawa, Godambe and Taylor (1997) contain extensive discussions on estimating equations. In this paper we introduce a generalized bootstrap technique for estimators obtained by solving estimating equations.

We use the following framework: Suppose $\{\phi_{ni}(Z_{ni},\beta),\ 1\leq i\leq n,\ n\geq 1\}$ is a triangular sequence of functions taking values in $\mathbb{R}^p$, $\{Z_{ni}\}$ being a sequence of observable random variables and $\beta\in\mathcal{B}\subset\mathbb{R}^p$. Assume that $E\phi_{ni}(Z_{ni},\beta_0)=0, 1\leq i\leq n,\ n\geq 1$ for some unique $\beta_0\in\mathcal{B}$. The "parameter" $\beta_0$ is unknown, and its estimator $\hat{\beta}_n$ is obtained by solving (often uniquely) the estimating equations

$$(1.1) \qquad \sum_{i=1}^{n}\phi_{ni}(Z_{ni},\beta)=0.$$

Typically, $\{\phi_{ni}(Z_{ni},\beta_0)\}$ form a triangular array of martingale differences.

The major objective of this paper is to estimate the sampling distribution and the asymptotic variance of $\hat{\beta}_n$ by a new approach to resampling. We define our resampling estimator $\hat{\beta}_B$ as the solution of

$$\sum_{i=1}^{n} w_{ni} \phi_{ni}(Z_{ni}, \beta) = 0, \tag{1.2}$$

where $\{w_{ni}, 1 \leq i \leq n, n \geq 1\}$ is a triangular sequence of random variables, independent of $\{Z_{ni}\}$. These are the "bootstrap weights." Note that essentially the same algorithm computes $\hat{\beta}_n$ and the Monte Carlo samples of $\hat{\beta}_B$. This makes the proposed bootstrap software friendly.

The normal equations $\sum \mathbf{x}_{ni}(y_{ni} - \mathbf{x}_{ni}^T \beta) = 0$ for the least squares estimator (LSE) in linear regression is a special case of (1.1). With $(w_{n1}, \ldots, w_{nn}) \sim$ Multinomial$(n, 1/n, \ldots, 1/n)$ we get the *paired bootstrap* (PB) estimator from (1.2). Other choices of $w_{ni}$'s yield the delete-$d$ jackknives, the Bayesian bootstrap, the $m$-out-of-$n$ bootstrap and variations of these. Hence we refer to resampling by (1.2) as the *generalized bootstrap* (GBS). Origins of the concept of resampling equations may be traced back to Freedman and Peters (1984) and Rao and Zhao (1992), where the bootstrap was carried out using equations, as distinguished from resampling observations or residuals. Note that the GBS technique is different from the bootstraps suggested by Lele (1991) and Hu and Kalbfleisch (2000) for estimating equations.

In Section 2.1 we state the conditions on GBS weights. In Section 2.2 we briefly discuss some examples of GBS schemes. Since every choice of distributions of the bootstrap weights corresponds to a different GBS technique, it is of interest to compare their relative performances. A theoretical comparison of different GBS techniques is under study, and some preliminary results may be found in Bose and Chatterjee (2002). Section 2.3 contains examples to illustrate the implementation of GBS. The standard GBS schemes, obtained by taking i.i.d. or multinomial weights, appear to perform competitively in a variety of problems, although there is some model and sample size dependent performance variation.

In Section 3.1 we assume $p = 1$ and establish asymptotic linearizations of $\hat{\beta}_n$ and $\hat{\beta}_B$. The distributional consistency of the GBS follows easily from these. In Section 3.3 we consider models with increasing dimension by letting $p \to \infty$ as $n \to \infty$ and establish similar results.

For the distribution of linear regression $M$-estimators, our results in this paper imply that the GBS is consistent even when regressors are random, errors are heteroscedastic or parameter dimension is increasing with sample size. This may be compared with Lahiri (1992), where a nonnaive residual bootstrap (RB) was found to be second-order accurate when covariates are nonrandom, errors are i.i.d. and the parameter dimension is fixed. While first-order consistency of GBS is achieved under relaxed assumptions, the GBS is second-order accurate only after a complicated bias correction and Studentization.

In Section 3.2 for dimension $p = 1$, we obtain an asymptotic representation for the GBS variance estimator, similar to the work of Liu and Singh (1992) and Hu (2001). Our result implies that for the asymptotic variance of linear regression $M$-estimators the GBS is consistent even when the errors are heteroscedastic, and yet can have greater asymptotic efficiency than some resampling schemes that are consistent only under homoscedasticity.

The technical framework used here is for estimating equations similar to $M$-estimation problems. However, the underlying principle of GBS may be applicable to a much wider class of statistical problems.

**2. GBS weights: conditions and examples.** In this section we spell out the technical conditions needed on the GBS weights and give examples of classes of weights which satisfy these conditions. We also illustrate the implementation of GBS through a few examples.

2.1. *Conditions on bootstrap weights.* Let $\{w_{ni};\ 1 \leq i \leq n,\ n \geq 1\}$ be a triangular array of nonnegative random variables such that for each $n$, the weights $w_{n1}, \ldots, w_{nn}$ are exchangeable. These are to be used as weights and *we drop the suffix $n$ from the notation.* $\mathbf{P}_B$ and $\mathbf{E}_B$, respectively, denote bootstrap probability and expectation conditional on the data. Let

$$V(w_i) = \sigma_n^2, \qquad W_i = (w_i - 1)/\sigma_n,$$

$$c_{ijk\cdots} = E(W_a^i W_b^j W_c^k \cdots) \quad \text{and} \quad a_n = \left[ \sup_{\|c\|=1} \sum_{i=1}^n E(c^T \phi_{ni})^2 \right]^{1/2}.$$

In the conditions below, $p$ is the dimension of the parameter space, which is allowed to tend to infinity with data size $n$ in Section 3.3.

The first set of conditions is fairly universal and is satisfied by all known examples of bootstrap weights.

BW (Basic conditions):

(2.1) $$Ew_i = 1,$$

(2.2) $$0 < \sigma_n^2 = o(\min(a_n^2 p^{-1}, n)),$$

(2.3) $$c_{11} = O(n^{-1}).$$

Schemes like the classical bootstrap and the delete-$d$ jackknife satisfy $\sum_{i=1}^n w_{ni} = C_n$ for some nonrandom sequence $\{C_n\}$. This implies that $c_{11} = -1/(n-1)$ and thus (2.3) is satisfied.

Additional assumptions required for distributional consistency are:

CLTW (Conditions for GBS CLT):

(2.4) $$c_{22} \to 1, \qquad c_4 < \infty.$$

For variance estimation, we need the basic conditions, (2.5), (2.6) and either part (a) or part (b) of (2.7) stated below.

Let $C^+ \subset (0, \infty)$ be a compact set, and let $\mathcal{W}$ be the set on which at least $m_0$ of the weights are greater than some fixed constant $k_2 > 0$.

VW (Conditions for GBS variance):

(2.5) $$\mathbf{P}_\mathrm{B}[\mathcal{W}] = 1 - O_P(n^{-1}),$$

(2.6) $$c_{i_1 i_2 \cdots i_k} = O(n^{-k+1} \sigma_n^{-1}) \qquad \forall\, i_1, i_2, \ldots, i_k \text{ satisfying } \sum_{j=1}^{k} i_j = 3,$$

(2.7)
(a) $\sigma_n^2 \in C^+$; $\quad c_{i_1 \cdots i_k} = O(\min(n^{-k+2}, 1)) \ \forall\, i_1, \ldots, i_k \text{ with } \sum_{j=1}^{k} i_j = 4,$

(b) $\sigma_n^2 \to 0$; $\quad c_{i_1 \cdots i_k} = O(n^{-k+2}) \ \forall\, i_1, \ldots, i_k \text{ with } \sum_{j=1}^{k} i_j = 4.$

In (2.6) and (2.7) the $i_j$'s are positive integers. In the following we refer to conditions (2.5), (2.6) and (2.7)(a) as VW(a) and to conditions (2.5), (2.6) and (2.7)(b) as VW(b).

2.2. *Examples of GBS weights.* We now list some common resampling techniques that are special cases of GBS.

(a) Suppose $\mathbf{w}_n = (w_{n1}, \ldots, w_{nn}) \sim \text{Multinomial}(n; 1/n, \ldots, 1/n)$. These weights can be interpreted as simple random sampling with replacement of the functionals to minimize and essentially correspond to the *classical bootstrap* of Efron (1979). Apart from BW, these weights also satisfy CLTW and VW(a).

Suppose instead that we select $m$ data points out of $n$ where typically $m \to \infty$ and $m/n \to 0$. If the selection is with replacement, the weights are an appropriately rescaled random sample from $\text{Multinomial}(m; 1/n, \ldots, 1/n)$. This scheme is usually called the *$m$-out-of-$n$ bootstrap*. If the selection is without replacement, the scheme can be identified with the *delete-$(n-m)$ jackknife*. For either situation, BW and CLTW hold.

See Præstgaard and Wellner (1993) for other variations and adaptations of the classical bootstrap.

(b) The *Bayesian bootstrap* [Rubin (1981)] and its variations [see Zheng and Tu (1988) and Lo (1991)] essentially use $\mathbf{w}_n \sim \text{Dirichlet}(\alpha, \ldots, \alpha)$. The *weighted likelihood bootstrap* of Newton and Raftery (1994) is also a variation, where $\phi_{ni}(\cdot)$ has a log-likelihood interpretation. The conditions BW, CLTW and VW(a) are satisfied.

(c) The *jackknives* are specially geared towards estimation of bias and variance. Suppose $\theta_n$ is an estimator based on $n$ observations and we wish to estimate its variance.

In its simplest form, the delete-1 jackknife estimator is obtained as follows: Drop the $i$th observation and recompute the estimator, say $\theta_{n,i}$, on the basis of the remaining $n-1$ observations. Then the jackknife estimator of the variance is $v = (n-1)n^{-1}\sum_{i=1}^{n}(\theta_{n,i} - \theta_n)^2$. To visualize the delete-1 jackknife as coming from a sequence of random weights, consider all vectors $\eta_i$, $1 \le i \le n$, of length $n$ where the $i$th coordinate of $\eta_i$ is zero and the rest are 1. Let $P(\mathbf{w}_n = n(n-1)^{-1}\eta_i) = 1/n$ for $1 \le i \le n$. The above estimator is then obtained after appropriate averaging over this uniform weight distribution.

The delete-$d$ jackknife deletes $d$ observations at a time and has a similar interpretation. If $n - d \to \infty$, then BW holds. If $d/n \to c \in (0,1)$, then CLTW and VW(a) hold. If $d/n \to 0$, then VW(b) holds.

The *downweight-$d$ jackknife* is a variation of the above. For $1 \le d \le n$ consider the $n$-dimensional vectors $\eta_{n:i_1,i_2,\ldots,i_d}$ where the $j$th coordinate of $\eta_{n:i_1,\ldots,i_d}$ is $d/n$ if $j$ is one of $i_1,\ldots,i_d$, else it is $(n+d)/n$. The resampling weights vector is a random sample from the set of $\eta$. The asymptotic properties of these weights are similar to the delete-$d$ jackknives. However, since no observation is assigned a weight zero, model assumptions like (3.16) are not needed.

2.3. *Examples on implementation of GBS in some models.* We consider three examples in this section. Important non-GBS techniques such as the RB and the wild bootstrap (WB) are also included for comparison.

EXAMPLE 2.1. *Heteroscedastic time series*: Consider the following model: $X_t = \phi X_{t-1} + e_t$, $t = 1, \ldots, n$, where $X_0 \equiv 0$, and $\{e_t\}$ is a sequence of independent, normal, mean-zero random variables with $Ee_t^2 = \sigma_1^2$ if $t$ is odd and $Ee_t^2 = \sigma_2^2$ if $t$ is even. Suppose that the unknown $\phi$ is estimated by the LSE $\hat{\phi} = \sum X_t X_{t-1} / \sum X_{t-1}^2$. Let $V_n = E(\sqrt{n}(\hat{\phi} - \phi))^2$ be the quantity of interest to be estimated using resampling techniques. In general $\phi$, $\sigma_1^2$ and $\sigma_2^2$ are unknown. For simulation purposes we let $\phi = 0.2$, $\sigma_1^2 = 1$ and $\sigma_2^2 = 100$.

We study the wild bootstrap (WB) [Wu (1986) and Mammen (1992)], GBS(1) with Multinomial$(n, 1/n, \ldots, 1/n)$ weights and GBS(2) with i.i.d. Uniform$(0.5, 1.5)$ weights. For simplicity, we use i.i.d. $N(0, 1)$ weights for the WB in all examples in this section. We used 10,000 simulations and bootstrap sample size of 1000 on each of the four resampling techniques.

In Table 1 the first column indicates the sample size. The value of $V_n$ depends on $n$, but is approximately 0.11 for all three $n$ values reported.

TABLE 1
*Mean (and variance) of estimates of $V_n$ from heteroscedastic AR(1) process (see Example 2.1) for residual bootstrap (RB), wild bootstrap (WB), GBS with multinomial weights [GBS(1)], GBS with Uniform(0.5, 1.5) weights [GBS(2)] over 10,000 simulation runs*

| $n$ | RB | WB | GBS(1) | GBS(2) |
|-----|------------|------------|------------|------------|
| 15  | 0.891 (0.013) | 0.151 (0.008) | 0.555 (0.466) | 0.126 (0.007) |
| 30  | 0.904 (0.006) | 0.146 (0.007) | 0.229 (0.021) | 0.131 (0.005) |
| 50  | 0.928 (0.004) | 0.124 (0.002) | 0.161 (0.005) | 0.121 (0.002) |

The first column denotes the sample size. Resample size is 1000.

The second column has the average over $k$ of $V_{\text{RB}}^k$, the residual bootstrap estimate of the variance of $\hat{\phi}$ for the $k$th simulation run. The variance of $V_{\text{RB}}^k$ is given in parentheses. The figures in columns three to five have similar interpretation for WB, GBS(1) and GBS(2).

From Table 1, it can be seen that RB, as expected, fails since it is not adapted for heteroscedasticity. GBS(1) is better, but is erratic at low sample sizes, a fact reflected in the high variance value of 0.47. WB does reasonably well, but is consistently outperformed by GBS(2). However, for larger sample sizes the difference between the latter three is nominal.

EXAMPLE 2.2. *Generalized linear models*: Suppose $\{Y_{ij},\ j=1,\ldots,N_i\}$ are independent Bernoulli($p_i(\beta)$) random variables with $p_i(\beta) = [1+\exp(t_i)]^{-1}\exp(t_i)$ and $t_i = \beta_0 + \beta_1 X_i$ for $i = 1,\ldots,n$.

We use $\{(N_i, X_i), i = 1,\ldots, n = 10\}$ from the data relating to effectiveness of ethylene oxide as a fumigant [Myers, Montgomery and Vining (2002), page 129]. Analysis of the actually observed $Y_i = \sum_{j=1}^{N_i} Y_{ij}$ values reported in those data yields maximum likelihood estimates $-17.90$ for $\beta_0$ and $6.28$ for $\beta_1$. We use these values as the true parameter values and simulate the $Y_{ij}$'s according to the model described above, and obtain estimates $\hat{\beta}_0$ and $\hat{\beta}_1$ of $\beta_0$ and $\beta_1$ by solving the likelihood equation.

We study the WB and two GBS techniques in this example. Let $N = \sum_{i=1}^{n} N_i$. We use $(w_1,\ldots,w_N)$ from Multinomial$(N; 1/N,\ldots, 1/N)$ for GBS(1) and $w_i$'s i.i.d. Exponential(1) for GBS(3). The WB is based on residuals, for which there are several choices. In the present example, we define $\hat{p}_{ij} = (Y_{ij} + \delta)/(1 + 2\delta)$ as the observed proportion of success. The constant $\delta = 0.001$ is used to avoid computational pathologies arising from the observed proportions being 0 or 1. Let $\tilde{t}_{ij} = \log(\hat{p}_{ij}/(1-\hat{p}_{ij}))$ be the "observed logit," while $\hat{t}_i = \hat{\beta}_0 + \hat{\beta}_1 X_i$ is the "fitted logit." Define $r_{ij} = \tilde{t}_{ij} - \hat{t}_i$ as the $ij$th residual, and $Y_{ij}^* = \hat{t}_i + U_{ij} r_{ij}$ where $U_{ij}$'s are i.i.d. $N(0,1)$ random variables.

Table 2
*Observed logit, average confidence interval length and coverage percentage from wild bootstrap (WB), GBS with multinomial weights [GBS(1)] and GBS with independent exponential weights [GBS(3)] for each of the 10 data points from Example 2.2*

| Case (Logit) | WB | GBS(1) | GBS(3) |
|---|---|---|---|
| 1 (2.264) | 0.96 (0.1) | 1.27 (96.4) | 1.25 (96.4) |
| 2 (2.213) | 0.94 (0.1) | 1.25 (96.4) | 1.23 (96.4) |
| 3 (1.791) | 0.84 (0.1) | 1.07 (96.3) | 1.05 (96.4) |
| 4 (1.220) | 0.70 (0.1) | 0.85 (95.2) | 0.83 (95.1) |
| 5 (1.099) | 0.67 (0.1) | 0.80 (95.2) | 0.79 (95.0) |
| 6 (0.321) | 0.61 (94.8) | 0.66 (92.8) | 0.66 (91.9) |
| 7 (−0.182) | 0.69 (99.7) | 0.75 (97.0) | 0.74 (97.2) |
| 8 (−0.567) | 0.76 (57.1) | 0.83 (95.8) | 0.83 (95.8) |
| 9 (−1.020) | 0.85 (9.6) | 0.97 (94.8) | 0.96 (94.9) |
| 10 (−2.956) | 1.33 (0.3) | 1.79 (95.2) | 1.76 (94.7) |

The nominal coverage is 95%. Resample size is 1000.

For each "true logit" $t_i = -17.90 + 6.28 X_i$, $i = 1, \ldots, n$, we obtain percentile based 95% confidence intervals using the three resampling schemes. Resampling size taken is $B = 1000$. This exercise is repeated $I = 1000$ times, and in Table 2 we report the average confidence interval length and coverage percentage over these 1000 replications of the experiment. The WB performs poorly in this example, since $r_{ij}$ depend on $\delta$ and carry little information on the variability of the data. The GBS techniques perform excellently in comparison.

Note that the likelihood is a function of the sufficient statistics $Y_i = \sum_{j=1}^{N_i} Y_{ij}$, and sometimes only the $Y_i$'s, and not the individual $Y_{ij}$'s, are available data. There we may use $(Y_i + \delta)/(N_i + 2\delta)$ as the "observed proportion of success" associated with the $i$th covariate. This improves the performance of the WB if $N_i$'s are large, and if $Y_i$'s are not close to zero or $N_i$. However, in many problems $N_i > 1$ may not be an available option.

EXAMPLE 2.3. *Nonlinear regression*: We consider the isomerization data from Huet, Bouvier, Gruet and Jolivet [(1996), page 11]. The reaction rate of the catalytic isomerization of $n$-pentane to isopentane depends on partial pressure at various stages. The model for the $i$th reaction rate $y_i$ is $y_i = f(X_i, \theta) + e_i$, where $f(X_i, \theta) = \frac{\theta_1 \theta_3 (P_i - I_i/1.632)}{1 + \theta_2 H_i + \theta_3 P_i + \theta_4 I_i}$ and $X_i = (H_i, P_i, I_i)^T$ are the corresponding partial pressure values. The $e_i$'s are i.i.d. random variables. The parameter $\theta = (\theta_1, \theta_2, \theta_3, \theta_4)^T$ is estimated by minimizing $\Psi_n(\theta) = \sum_{i=1}^n (y_i - f(X_i, \theta))^2$ with the resulting estimate $\hat{\theta} = (35.9193, 0.0708583, 0.0377385, 0.167166)^T$.

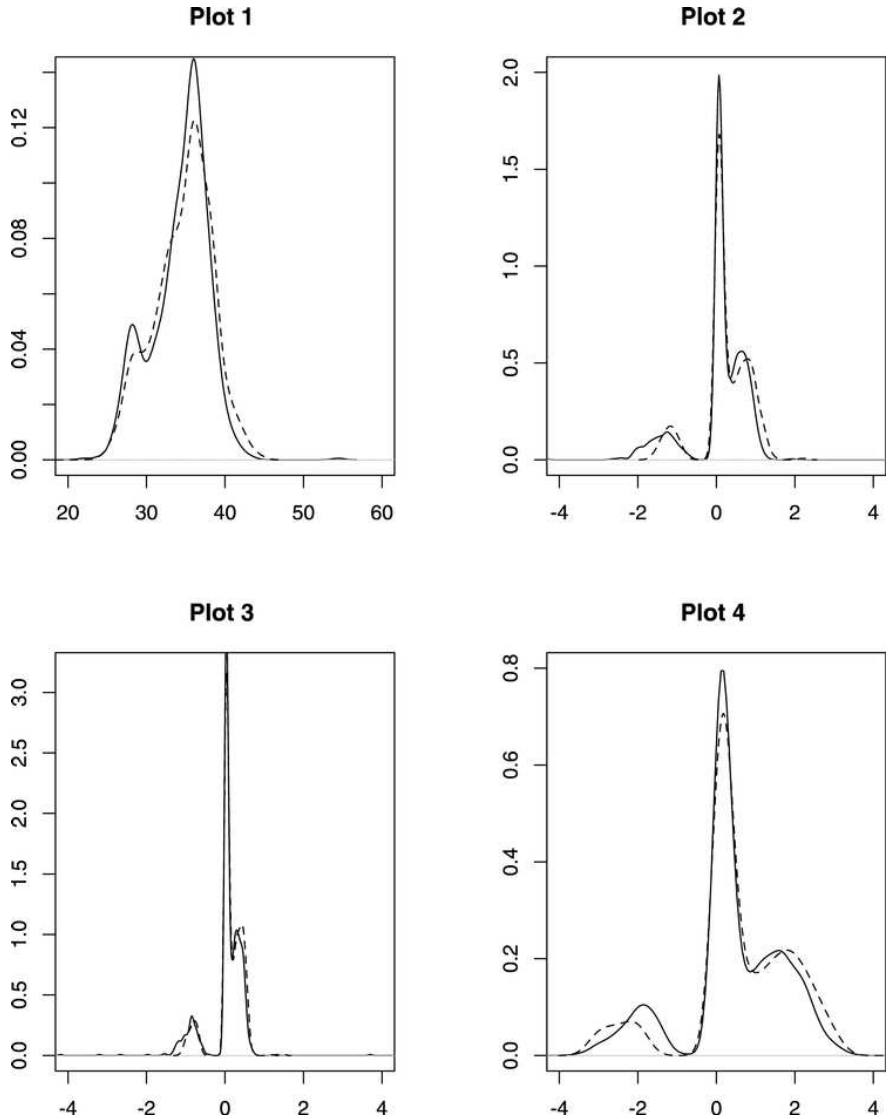

FIG. 1. *Plots of GBS with multinomial weights* [GBS(1)] *(solid line) and residual bootstrap (RB) (broken line) densities for the four parameters in Example* 2.3. *Plot i corresponds to* $\theta_i$, $i = 1, 2, 3, 4$. *Resample size is* 1000.

The analysis in Huet, Bouvier, Gruet and Jolivet (1996) includes an RB using Studentized quantities for each $\theta_i$, and the resulting 95% equal-tail confidence interval does not include zero for any of the $\theta_i$'s.

We study the RB, WB, GBS(1) and GBS(3) here. Note that $(w_1, \ldots, w_n) \sim$ Multinomial$(n; 1/n, \ldots, 1/n)$ for GBS(1) and the $w_i$'s are i.i.d. Exponential(1) for GBS(3). Figure 1 represents the density histograms from RB and GBS(1) overlaid on each other. Notice that for each $\theta_i$ the resampling densities have two prominent modes, one near the estimate $\hat{\theta}$ and the other near $\theta^* = (33.343956, -1.84281206, -1.0338937, -4.31406116)^T$. Note that $\Psi_n(\theta^*) = 3.26$, a value quite close to $\Psi_n(\hat{\theta}) = 3.23$. The results from GBS(3) are similar to those of GBS(1), while in WB the peak at $\theta^*$ is slightly less prominent.

The estimates $\hat{\theta}$ and $\theta^*$ represent two substantially different chemical processes. This being real data, it is not known if $\hat{\theta}$ or $\theta^*$ is closer to $\theta$. However, the presence of $\theta^*$ is not revealed in the analysis of Huet et al. The bimodal curves in Figure 1 suggest that convex confidence intervals make a bad summarization in the present problem. A less sensitive bootstrap such as GBS may thus be useful in revealing features in data that theoretically superior but sensitive resampling techniques may miss.

**3. Main results.** In Section 3.1 we assume $p = 1$ and obtain asymptotic representations of $\hat{\beta}_n$ and $\hat{\beta}_B$ in Theorems 3.1 and 3.2. This establishes the consistency of the GBS for estimating the distribution. In Section 3.3 we consider general $p$, including the case where $p \to \infty$ as $n \to \infty$ and obtain similar results in Theorems 3.4 and 3.5.

In Section 3.2 we focus on the variance estimation problem. We assume that the $\phi_{ni}, 1 \leq i \leq n$, are independent and $p = 1$. In Theorem 3.3 we establish an asymptotic representation of the GBS variance estimator, thereby generalizing part of the work of Liu and Singh (1992) and Hu (2001). All proofs are only sketched and complete details are available from the authors.

We discuss specific model conditions in the respective Sections. We introduce some of the notation here: throughout, $k$ and $K$, with or without suffix, are used as generic constants. Two conventions are used: any condition stated for a random function is assumed to hold almost surely unless otherwise stated; and "for all $\beta$" always means for all $\beta$ in an open neighborhood of $\beta_0$.

Write $\phi_{ni}(Z_{ni}, \beta) = \phi_{ni}(\beta) = (\phi_{ni(1)}(\beta), \ldots, \phi_{ni(p)}(\beta))^T$. Thus the $a$th coordinate of $\phi_{ni}$ is $\phi_{ni(a)}$, $a = 1, \ldots, p$. Let $\phi_{0ni}(\beta) \equiv \phi_{ni}(\beta)$ and for $k \geq 0$, $\phi_{(k+1)ni(a)}(\beta) = \frac{\partial}{\partial \beta} \phi_{kni(a)}(\beta)$. Let $\phi_{kni(a)} = \phi_{kni(a)}(\beta_0)$. For each $\phi_{ni(a)}(\beta)$, we assume that the following Taylor series expansion holds:

$$(3.1) \qquad \phi_{ni(a)}(\beta + t) = \phi_{ni(a)}(\beta) + \phi_{1ni(a)}^T(\beta)t + 2^{-1}t^T H_{2ni(a)}(\beta_1)t$$

for $\beta_1 = \beta + ct$ and for some $0 < c < 1$.

3.1. *Asymptotics and bootstrap for $p = 1$.* When $p = 1$, we simplify notation by suppressing the last index and thus: $\phi_{ni}(Z_{ni}, \beta) \equiv \phi_{ni(1)}(\beta)$, $\phi_{0ni}(\beta) \equiv \phi_{ni}(\beta)$, $\phi_{(k+1)ni}(\beta) = \frac{\partial}{\partial \beta} \phi_{kni}(\beta)$ and $\phi_{kni} = \phi_{kni}(\beta_0)$. We then write (3.1) as $\phi_{ni}(\beta + t) = \phi_{ni}(\beta) + \phi_{1ni}(\beta)t + 2^{-1}\phi_{2ni}(\beta_1)t^2$ for $\beta_1 = \beta + ct$ and for some $0 < c < 1$.

Let

$$\gamma_{1n}^2 = a_n^{-2} \sum_{i=1}^n E\phi_{1ni}, \qquad S_{nj} = \sum_{i=1}^j \phi_{ni}.$$

Assume that $a_n = [\sum_{i=1}^n E\phi_{ni}^2]^{1/2} \to \infty$.

*Assumptions for Section* 3.1. Assume that for every $n$, there is a sequence of $\sigma$-fields $\mathcal{F}_{n1} \subset \cdots \mathcal{F}_{nn}$, such that $\{S_{nj}, \mathcal{F}_{nj}, 1 \leq j \leq n\}$ is a martingale. Further, with $\eta_n = \max(\sigma_n^2, 1)$,

(3.2) $\qquad\qquad\qquad E\phi_{ni} = 0 \qquad$ for all $1 \leq i \leq n, \quad n \geq 1,$

(3.3) $\qquad\qquad\qquad 0 < k_2 < \gamma_{1n},$

(3.4) $\qquad E\Big[\sum(\phi_{1ni} - E\phi_{1ni})\Big]^2 = o(a_n^4 \eta_n^{-1}).$

There exist $\delta_0 > 0$ and $M_{2ni}$ such that

(3.5) $\qquad \sup_{|\beta - \beta_0| < \delta_0} |\phi_{2ni}(\beta)| \leq M_{2ni} \quad$ and $\quad E\Big(\sum_{i=1}^n M_{2ni}\Big)^2 = o(a_n^6 \eta_n^{-1}).$

The triangular sequence $X_{ni} = (\sum_{i=1}^n E\phi_{ni}^2)^{-1/2} \phi_{ni}$ satisfies

(3.6) $\qquad\qquad \sum_{i=1}^n X_{ni}^2 \xrightarrow{p} 1, \qquad E\Big(\max_i |X_{ni}|\Big) \to 0.$

THEOREM 3.1. *Under* (3.2)–(3.5) *there exists a sequence* $\{\hat{\beta}_n\}$ *of solutions of* (1.1) *such that*

(3.7) $\qquad\qquad a_n(\hat{\beta}_n - \beta_0) = O_P(1),$

(3.8) $\qquad\qquad a_n\gamma_{1n}(\hat{\beta}_n - \beta_0) = -a_n^{-1}\sum_{i=1}^n \phi_{ni} + r_n,$

*where* $r_n = o_P(1)$. *In addition, if* (3.6) *holds, then* $[\sum_{i=1}^n E\phi_{1ni}]^{1/2}(\hat{\beta}_n - \beta_0) \xrightarrow{\mathcal{D}} N(0,1)$.

Before we give the proof, we note that in general a sequence of solutions need not be measurable. See, for example, Ferguson (1996). However, there are enough assumptions in our model to guarantee this measurability. We omit these arguments here and also for the subsequent results.

PROOF OF THEOREM 3.1. Fix any $\varepsilon > 0$. By Chebyshev's inequality and (3.2), there exists a $K > 0$ such that

$$\text{Prob}\left[\left|a_n^{-1}\sum_{i=1}^n \phi_{ni}\right| > K\right] < \varepsilon/2. \tag{3.9}$$

Define $S_n(t) = a_n^{-1}\sum_{i=1}^n[\phi_{ni}(\beta_0 + a_n^{-1}t) - \phi_{ni}(\beta_0)] - \gamma_{1n}t$. Using a Taylor series expansion of $\phi_{ni}(\cdot)$ about $\beta_0$ and (3.4)–(3.5), we can show that, given any constant $C > 0$, for all large $n$,

$$E\left[\sup_{|t|\leq C}|S_n(t)|\right]^2 = o(1). \tag{3.10}$$

Now note that

$$\inf_{|t|=C}\left\{a_n^{-1}t\sum_{i=1}^n\phi_{ni}(\beta_0 + a_n^{-1}t)\right\}$$

$$\geq -C\sup_{|t|=C}|S_n(t)| + C^2\gamma_{1n} - Ca^{-1}\left|\sum_{i=1}^n\phi_{ni}\right|. \tag{3.11}$$

From (3.9)–(3.11) we have, choosing $C$ large enough,

$$\text{Prob}\left[\inf_{|t|=C}\left\{a_n^{-1}t\sum_{i=1}^n\phi_{ni}(\beta_0 + a_n^{-1}t)\right\} > 0\right]$$

$$\geq \text{Prob}\left[a_n^{-1}\left|\sum_{i=1}^n\phi_{ni}\right| + \sup_{|t|=C}|S_n(t)| \leq C\gamma_{1n}\right]$$

$$= 1 - \text{Prob}\left[a_n^{-1}\left|\sum_{i=1}^n\phi_{ni}\right| + \sup_{|t|=C}|S_n(t)| > C\gamma_{1n}\right]$$

$$\geq 1 - \text{Prob}\left[a_n^{-1}\left|\sum_{i=1}^n\phi_{ni}\right| > Ck_2/4\right] - \text{Prob}\left[\sup_{|t|=C}|S_n(t)| > Ck_2/4\right]$$

$$\geq 1 - \varepsilon \quad \text{for all } n \text{ sufficiently large.}$$

By the continuity of $\sum_{i=1}^n \phi_{ni}(\beta)$ in $\beta$, this means that, for fixed $\varepsilon > 0$ for all $n$ sufficiently large, there exists a $C$ such that

$$\sum_{i=1}^n \phi_{ni}(\beta_0 + a_n^{-1}t) = 0 \quad \text{has a root } t = T_n \text{ in } |t| \leq C \text{ with probability } > 1 - \varepsilon.$$

Defining $\hat{\beta}_n = \beta_0 + a_n^{-1}T_n$ when such $T_n$ exists and as an arbitrary zero of $\sum_{i=1}^n \phi_{ni}(\beta) = 0$ otherwise, we get a solution to (1.1) which satisfies, for fixed $\varepsilon > 0$, $\text{Prob}[a_n|\hat{\beta}_n - \beta_0| \leq C] \geq 1 - \varepsilon$ for all $n$ large enough. This shows (3.7).

Now with this $C$ fixed, by arguments similar to those of (3.10), we obtain that $a_n\gamma_{1n}(\hat{\beta}_n - \beta_0) = -a_n^{-1}\sum_{i=1}^n \phi_{ni} + r_n$, where $r_n = o_P(1)$.

The dependence of $\hat{\beta}_n$ on the choice of $\varepsilon$ may be taken care of as described in Serfling [(1980), page 148]. Briefly, this is as follows:

Since $\hat{\beta}_n \equiv \hat{\beta}_{n,\varepsilon}(\omega) \xrightarrow{p} \beta_0$ there is a subsequence along which the convergence is with probability 1, and we may restrict attention to this subsequence only. Thus in our definition of $\hat{\beta}_{n,\varepsilon}(\omega)$ above, for every $\varepsilon > 0$ there is an $N_\varepsilon$ such that for all $n > N_\varepsilon$, $\omega$ belongs to a probability-1 set $\Omega_\varepsilon$. Then on the probability-1 set $\Omega_0 = \bigcap_{k \geq 1} \Omega_{1/k}$, without loss of generality we have a nondecreasing sequence of integers $N_1(\omega) \leq N_{1/2}(\omega) \leq \cdots \leq N_{1/k}(\omega)\ldots$. For $n \in [N_{1/k}(\omega), N_{1/(k+1)}(\omega))$, we define $\hat{\beta}_n = \hat{\beta}_{n,1/k}(\omega)$ and let $\hat{\beta}_n = 0$ otherwise. Then the new sequence $\{\hat{\beta}_n\}$ has all the desired properties.

Further, assumption (3.6) ensures that $\sum_i X_{ni} = a_n^{-1} s_n^{-1} \sum_i \phi_{ni} \xrightarrow{\mathcal{D}} N(0,1)$ by Theorem 5.4.2 of Borovskikh and Korolyuk (1997). □

Henceforth we work with that sequence of solutions $\{\hat{\beta}_n\}$ which satisfies Theorem 3.1.

The bootstrap estimator is obtained by solving (1.2). The next theorem is on its asymptotic representation and consistency. Let

$$F_n(x) = P\left[\left[\sum_{i=1}^n E\phi_{1ni}\right]^{1/2} (\hat{\beta}_n - \beta_0) \leq x\right],$$

$$F_{Bn}(x) = \mathbf{P}_{\mathrm{B}}\left[\sigma_n^{-1}\left(\sum_{i=1}^n \phi_{1ni}(\hat{\beta}_n)\right)^{1/2} (\hat{\beta}_B - \hat{\beta}_n) \leq x\right].$$

THEOREM 3.2. *Assume* (3.2)–(3.5) *and that the bootstrap weights satisfy BW. Then there exists a sequence* $\{\hat{\beta}_B\}$ *of solutions of* (1.2) *such that*

$$(3.12)\quad \sigma_n^{-1}\left[\sum_{i=1}^n \phi_{1ni}(\hat{\beta}_n)\right]^{1/2} (\hat{\beta}_B - \hat{\beta}_n) = -a_n^{-1}\sum_{i=1}^n W_i \phi_{ni}(\hat{\beta}_n)\phi_{ni} + r_{nB},$$

*where* $\mathbf{P}_{\mathrm{B}}(|r_{nB}| > \varepsilon) = o_P(1)$ *for any* $\varepsilon > 0$. *If in addition* (3.6) *and CLTW hold, then*

$$(3.13)\qquad \sup_x |F_{Bn}(x) - F_n(x)| \to 0 \quad \text{in probability.}$$

PROOF. The technique used in proving (3.12) is similar to the proof of (3.7) and (3.8), and we omit some of the details here.

Define $\hat{\gamma}_{1n} = a_n^{-2}\sum_{i=1}^n \phi_{1ni}(\hat{\beta}_n)$ and $S_{nB}(t) = a_n^{-1}\sum_{i=1}^n w_i[\phi_{ni}(\hat{\beta}_n + a_n^{-1}t) - \phi_{ni}(\hat{\beta}_n)] - \hat{\gamma}_{1n}t$. By arguments similar to those in the proof of Theorem 3.1

we have

$$\mathbf{P}_{\mathrm{B}}\left[\inf_{|t|=C\sigma_n}\left\{a_n^{-1}t\sum_{i=1}^n w_i\phi_{ni}(\hat{\beta}_n+a_n^{-1}t)\right\}>0\right]$$
$$\geq 1-\mathbf{P}_{\mathrm{B}}\left[\left|a_n^{-1}\sum_{i=1}^n w_i\phi_{ni}(\hat{\beta}_n)\right|>C\hat{\gamma}_{1n}\sigma_n/2\right]$$
$$-\mathbf{P}_{\mathrm{B}}\left[\inf_{|t|=C\sigma_n}|S_{nB}(t)|>C\hat{\gamma}_{1n}\sigma_n/2\right]$$
$$=1-U_{1C}-U_{2C} \quad \text{say.}$$

For given $\varepsilon>0$ and $\delta>0$, one can fix $C$ large enough such that for all $n$ sufficiently large, we have $\mathrm{Prob}[U_{iC}>\varepsilon/2]<\delta/2, i=1,2$.

Then with some algebra it can be established that $\sup_{|t|\leq C\sigma_n}|S_{nB}(t)|=\sigma_n r_{nB}$, where $\mathbf{P}_{\mathrm{B}}(|r_{nB}|>\varepsilon)=o_P(1)$ for any $\varepsilon>0$; then it follows that $\hat{\gamma}_{1n}\sigma_n^{-1}a_n(\hat{\beta}_B-\hat{\beta}_n)=-a_n^{-1}\sum_{i=1}^n W_i\phi_{ni}(\hat{\beta}_n)+r_{nB}$. This completes the proof of the first part. The second part follows from Theorem 3.1, the first part, and Lemma 4.6 of Præstgaard and Wellner (1993). We omit the details. □

3.2. *Asymptotics of the bootstrap variance estimator.* The estimation of the asymptotic variance of $\hat{\beta}_n$ is an important practical problem. In general, distributional convergence and variance estimation are different problems. For example, the delete-$d$ jackknife ($d/n \to 0$) is not distributionally consistent but is variance consistent for the i.i.d. sample mean. In this section we establish consistency of the GBS variance estimator via an asymptotic representation.

*Assumptions for Section* 3.2. We assume that the parameter is real valued ($p=1$), and that the $\{\phi_{ni}\}$ are independent. Also assume that

$$(3.14) \quad \phi_{ni}(\beta+t)=\phi_{ni}(\beta)+\phi_{1ni}(\beta)t+2^{-1}\phi_{2ni}(\beta)t^2+R_{ni}(t,\beta)t^2,$$

where $|R_{ni}(t,\beta)|<k|t|^\alpha$ for each $\beta$ for some $0<\alpha\leq 1$.

Assume that with $L=8(1+\alpha)$:

$$(3.15) \qquad \sum_{i=1}^n E|\phi_{ni}|^L + \sum_{i=1}^n E|\phi_{1ni}|^L + \sum_{i=1}^n E|\phi_{2ni}|^L = O(n).$$

Suppose $m_0$ is a specified integer, related to assumption (2.5) on the bootstrap weights. For any integer $m$ in $\{m_0,\ldots,n\}$ consider the subset $\mathcal{I}_m=\{i_1,i_2,\ldots,i_m\}$ of $\{1,2,\ldots,n\}$. We assume

$$(3.16) \qquad m^{-1}\sum_{i\in\mathcal{I}_m}\phi_{1ni}(\beta)>k_1>0$$

for every such choice of subset $\mathcal{I}_m$ of size $m$ from $\{1, 2, \ldots, n\}$, for every $m$ in $[m_0, n]$ and $\beta$ satisfying $\|\beta - \beta_0\| < \delta$ for a $\delta > 0$.

Resampling schemes like the PB and the delete-$d$ jackknives effectively select subsets of the data in the resample, and the model assumption (3.16) is required to hold on these subsets to make such resampling schemes feasible. See Wu (1986) and its discussion for more details on this. Assumption (3.16) helps in showing that under appropriate conditions the probability of a "bad" subset selection by the bootstrap or jackknife mechanism is small; see Proposition 3.1 (proof omitted). Some bootstrap clone methods and the downweight-$d$ jackknives do not require assumption (3.16). The assumptions above are not the most general for consistency. However, the stronger assumptions allow for more transparent computations.

PROPOSITION 3.1. *Assume the $\phi_{ni}$ are independent satisfying* (3.14) *with* (3.2)–(3.5), (3.15) *and* (3.16). *Assume $\hat{\beta}_n$ is a solution to* (1.1) *from Theorem* 3.1. *Let $\mathcal{A}$ be the set on which $m^{-1} \sum_{i \in \mathcal{I}_m} \phi_{1ni}(\hat{\beta}_n) > k_1/2 > 0$ for every such choice of subset $\mathcal{I}_m$ of size $m$ from $\{1, 2, \ldots, n\}$ and for every $m$ in $[m_0, n]$. Then* $\mathrm{Prob}[\mathcal{A}] > 1 - O(n^{-2})$.

For this section we define our bootstrap estimator $\hat{\beta}_B$ to be the solution to (1.2) on the set $\mathcal{A} \cap \mathcal{W}$, and $\hat{\beta}_n$ otherwise. This is to facilitate variance computations, and the minor alteration in the definition is of negligible consequence in the asymptotics. The set $\mathcal{A}$ is defined in Proposition 3.1, and $\mathcal{W}$ is defined in Section 2. The GBS variance estimate is $\mathbf{V}_{\mathrm{GBS}} = \sigma_n^{-2} \mathbf{E}_B (\hat{\beta}_B - \hat{\beta}_n)^2$. Note that the asymptotic variance of $n^{1/2} g_{1n} (\hat{\beta}_n - \beta_0)^2$ is $v_n = n^{-1} \sum_{i=1}^n E \phi_{ni}^2$. In the statement of the next theorem we have used $\phi$, $\phi_1$, $\phi_2$, respectively, for $\phi_{ni}$, $\phi_{1ni}$, $\phi_{2ni}$. The sums range from 1 to $n$. Also let $g_{1n} = n^{-1} \sum E \phi_1$, $g_{2n} = n^{-1} \sum E \phi_2$.

THEOREM 3.3. *Assume the $\phi_{ni}$ are independent satisfying* (3.14) *with* (3.2)–(3.5), (3.15) *and* (3.16). *Assume $\hat{\beta}_n$ is a solution to* (1.1) *from Theorem* 3.1. *Suppose the weights satisfy BW and either* VW(a) *or* VW(b). *Then*

$$
\begin{aligned}
ng_{1n}^2 (\mathbf{V}_{\mathrm{GBS}} - v_n) \\
= n^{-1} \sum (\phi^2 - E\phi^2) - \frac{2}{n^2 g_{1n}} \sum \phi \sum \phi \phi_1 \\
- \frac{2}{n^2 g_{1n}} \sum \phi^2 \sum (\phi_1 - E\phi_1) \\
+ \frac{2}{n^3 g_{1n}^2} \sum \phi \sum \phi^2 \sum \phi_2 + O_P(n^{-1}).
\end{aligned}
$$
(3.17)

The terms on the right-hand side of (3.17) are $O_P(n^{-1/2})$, so Theorem 3.3 shows in particular that the resampling variance of $n^{1/2}\sigma_n^{-1}(\hat{\beta}_B - \hat{\beta}_n)$ is consistent for the asymptotic variance of $n^{1/2}(\hat{\beta}_n - \beta)$.

REMARK 1. The above asymptotic representation is actually that of the mean squared error. However, the bias is of a negligible order compared to the variance, and thus the same representation holds for the asymptotic variance.

REMARK 2. For the least squares estimator in linear regression, $\phi_1$ is a constant and consequently $\phi_2$ is zero. There, using expansions for resampling variances, Liu and Singh (1992) classified resampling techniques in two groups: some are consistent even if errors are heteroscedastic, thus they are "robust" ($R$-class); others work only under homoscedasticity but have greater "efficiency" ($E$-class) than $R$-class techniques. Later, Bose and Kushary (1996) and Hu (2001) showed that the above classification breaks down if some other $M$-estimators are used.

Representation (3.17) is the same [up to $O_P(n^{-1})$ terms] as the $R$-class representation obtained for the PB for LSE in Liu and Singh [(1992), Theorem 2(ii)] and for general regression $M$-estimators in Hu [(2001), Theorem 2.2(ii)]. Note, however, representation (3.17) holds for a much broader class of problems than regression $M$-estimators.

By computations similar to those in Hu (2001), it can be shown that for particular choices of $\psi(\cdot)$ the GBS can be simultaneously robust against heteroscedasticity of errors as well as more efficient than $E$-class techniques.

PROOF OF THEOREM 3.3. We omit some of the details of the algebra involved in this proof. They are similar to those of Theorems 3.1 and 3.2.

Let us concentrate on the set $\mathcal{A} \cap \mathcal{W}$ only, since the contribution from the complement of this set is negligible. Define

$$U_{nB}(t) = \sigma_n^{-1} n^{-1/2} \sum_{i=1}^n w_i [\phi_{ni}(\hat{\beta}_n + \sigma_n n^{-1/2} t) - \phi_{ni}(\hat{\beta}_n)]$$

$$- n^{-1} t \sum_{i=1}^n w_i \phi_{1ni}(\hat{\beta}_n) - 2^{-1} \sigma_n n^{-3/2} t^2 \sum_{i=1}^n w_i \phi_{2ni}(\hat{\beta}_n).$$

Working along the lines of the proof of Theorem 3.2, we can show that

$$\mathbf{E}_B \left[ \sup_{|t| \leq C\sigma_n} |U_{nB}(t)| \right]^2 = O_P(n^{-(1+\alpha)}).$$

Now, under $\mathcal{A} \cap \mathcal{W}$ we may plug in $t = \sigma_n^{-1} n^{1/2}(\hat{\beta}_B - \hat{\beta}_n)$ in $U_{nB}(t)$, and after quite a lot of algebra we arrive at

$$-g_{1n}\sigma_n^{-1} n^{1/2}(\hat{\beta}_B - \hat{\beta}_n)$$
$$= n^{-1/2}\sum W_i \phi_{ni} - n^{-2} g_{1n}^{-1} \sum \phi_{ni} \sum W_i \phi_{1ni}$$
$$\quad - n^{-2} g_{1n}^{-1} \sum (\phi_{1ni} - E\phi_{1ni}) \sum W_i \phi_{ni}$$
$$\quad + n^{-5/2} g_{1n}^{-2} \sum \phi_{ni} \sum \phi_{2ni} \sum W_i \phi_{ni}$$
$$\quad + \sigma_n n^{-3/2} g_{1n}^{-1} \sum W_i \phi_{ni} \sum W_i \phi_{1ni}$$
$$\quad + 2^{-1} \sigma_n n^{-3/2} g_{2n} g_{1n}^{-2} \Big(\sum W_i \phi_{ni}\Big)^2 + R_{nB}$$
$$= C_n + T_{1n} + T_{2n} + T_{3n} + T_{4n} + T_{5n} + R_{nB} \qquad \text{say,}$$

where $\mathbf{E}_B R_{nB}^2 = O_P(n^{-(1+\alpha)})$.

Now it can be easily checked that $\mathbf{E}_B C_n^2 = O_P(1)$, and $\mathbf{E}_B T_{in}^2 = O_P(n^{-1})$, for $i = 1, \ldots, 5$. In the cross product, by direct computation $\mathbf{E}_B C_n T_{in} = O_P(n^{-1})$ for $i = 4, 5$, and hence $n g_{1n}^2 \mathbf{V}_{\text{GBS}} = \mathbf{E}_B C_n^2 + 2\mathbf{E}_B C_n(T_{1n} + T_{2n} + T_{3n}) + O_P(n^{-1})$. The rest of the proof follows by calculating the above moments. □

3.3. *Dimension asymptotics.* In this section we generalize the results of Section 3.1 to dimensions greater than 1 and also allow dimension $p = p_n \to \infty$ as the data size $n \to \infty$. Dimension asymptotics has been a major aspect of the study of resampling in the framework of linear regression [Bickel and Freedman (1983) and Mammen (1989, 1993)]. The classical residual-based bootstrap has been studied for the LSE [Bickel and Freedman (1983)] and for general $M$-estimators [Mammen (1989)] using nonrandom design matrices. The random design case and resampling using PB and WB have been studied in Mammen (1993). This section is an attempt to explore the high-dimensionality aspect in more general problems.

*Assumptions for Section* 3.3. The following notation will be used: $\|c\|$ is the Euclidean norm of a vector $c$, $A^T$ is the transpose of the matrix $A$, $\lambda_{\max}(A)$ and $\lambda_{\min}(A)$ are, respectively, the maximum and minimum eigenvalue of $A$.

Assume that

$$a_n = \Bigg[\sup_{\|c\|=1} \sum_{i=1}^n E(c^T \phi_{ni})^2\Bigg]^{1/2} \to \infty \qquad \text{as } n \to \infty.$$

Let
$$S_{nj} = \sum_{i=1}^{j} \phi_{ni} \quad \text{and} \quad \gamma_{1n}^2 = a_n^{-2} \sum_{i=1}^{n} E\phi_{1ni}.$$

Assume that for every $n$ there is a sequence of $\sigma$-fields $\mathcal{F}_{n1} \subset \cdots \subset \mathcal{F}_{nn}$, such that $\{S_{nj}, \mathcal{F}_{nj}, j = 1, \ldots, n\}$ is a martingale sequence. Recall that $\eta_n = \max(\sigma_n^2, 1)$.

Assume that

(3.18) $$E\phi_{ni} = 0,$$

(3.19) $$\sum_{i=1}^{n} \sum_{a=1}^{p} E\|\phi_{1ni(a)} - E\phi_{1ni(a)}\|^2 = o(a_n^4 \eta_n^{-1}).$$

For the symmetric matrix $H_{2ni(a)}$ in (3.1), for some $\delta_0 > 0$ there exists a symmetric matrix $M_{2ni(a)}$ such that

(3.20) $$\sup_{\{t:\|t\| \leq \delta_0\}} H_{2ni(a)}(\beta_0 + t) < M_{2ni(a)},$$

(3.21) $$\sum_{i=1}^{n} \sum_{a=1}^{p} E\lambda_{\max}^2(M_{2ni(a)}) = o(a_n^6 p^{-1} n^{-1} \eta_n^{-1}).$$

Let $\phi_{1ni}(\beta)$ be the $(p \times p)$ matrix, whose $a$th row is given by $\phi_{1ni(a)}^T(\beta)$, for $a = 1, \ldots, p$. Let $\Gamma_{1n}(\beta) = a_n^{-2} \sum_{i=1}^{n} \phi_{1ni}(\beta)$. Let $G_{1n} = a_n^{-2} \sum_{i=1}^{n} E\phi_{1ni}$. Assume

(3.22) $$0 < k_2 < \lambda_{\min}(G_{1n}).$$

Let $\{c = c_n \in \mathbb{R}^{p_n} = \mathbb{R}^p, \|c\| = 1\}$ be a fixed sequence of vectors on the unit balls of $p = p_n$-dimensional Euclidean spaces. Let

$$s_n^2 = p^{-2}\left[\left(\sum_{i=1}^{n} E\phi_{1ni}\right)^{-1} c\right]^T \left[\sum_{i=1}^{n} E\phi_{ni}\phi_{ni}^T\right]\left[\left(\sum_{i=1}^{n} E\phi_{1ni}\right)^{-1} c\right],$$

$$\hat{s}_n^2 = p^{-2}\left[\left(\sum_{i=1}^{n} \phi_{1ni}(\hat{\beta}_n)\right)^{-1} c\right]^T \left[\sum_{i=1}^{n} \phi_{ni}(\hat{\beta}_n)\phi_{ni}^T(\hat{\beta}_n)\right]\left[\left(\sum_{i=1}^{n} \phi_{1ni}(\hat{\beta}_n)\right)^{-1} c\right],$$

$$X_{ni} = -s_n^{-1}\left[\left(\sum_{i=1}^{n} E\phi_{1ni}\right)^{-1} c\right]^T \phi_{ni}.$$

Then $X_{ni}$ is measurable with respect to $\mathcal{F}_{ni}$ and satisfies

(3.23) $$\sum_{i=1}^{n} X_{ni}^2 \xrightarrow{p} 1, \qquad E\left(\max_i |X_{ni}|\right) \to 0.$$

THEOREM 3.4. *Under* (3.18)–(3.22) *there exists a sequence* $\{\hat{\beta}_n\}$ *of solutions of* (1.1) *such that if* $pa_n^{-2} \to 0$, *then*

$$(3.24) \qquad a_n p^{-1/2}\|(\hat{\beta}_n - \beta_0)\| = O_P(1),$$

$$(3.25) \qquad a_n p^{-1/2} G_{1n}^T(\hat{\beta}_n - \beta_0) = -a_n^{-1} p^{-1/2} \sum_{i=1}^n \phi_{ni} + r_n,$$

*where* $\|r_n\| = o_P(1)$. *In addition, if* (3.23) *is satisfied, then* $s_n^{-1} c^T(\hat{\beta}_n - \beta_0) \stackrel{\mathcal{D}}{\Rightarrow} N(0,1)$.

REMARK 3. The conditions (3.18)–(3.22) are nearly the same as conditions (C.1)–(C.3) of Lahiri (1992) except that he requires finite third moments for deriving Edgeworth expansions. It may also be noted that in most applications the $\phi_{1ni}$'s are uniformly almost surely bounded away from zero. Thus condition (3.22) [and (3.3)] is easily satisfied.

PROOF OF THEOREM 3.4. We first establish that given any $\varepsilon > 0$, $\exists K > 0$ such that $\text{Prob}[\|a_n^{-1} p^{-1/2} \sum_{i=1}^n \phi_{ni}\| > K] < \varepsilon/2$ using Chebyshev's inequality, (3.18) and that $\sum_{i=1}^n E(\|\phi_{ni}\|^2) = O(a_n^2 p)$. Let

$$S_n(t) = a_n^{-1} p^{-1/2} \sum_{i=1}^n [\phi_{ni}(\beta_0 + a_n^{-1} p^{1/2} t) - \phi_{ni}(\beta_0)] - G_{1n}^T t,$$

$$M_{1n} = \sum_{i,j=1}^n \sum_{a=1}^p (\phi_{1nia} - E\phi_{1nia})(\phi_{1nja} - E\phi_{1nja})^T.$$

Since $p/a_n^2 \to 0$, for every fixed $t$ eventually $\beta_0 + a_n^{-1} p^{1/2} t$ lies in the set $\{\beta_0 + x : \|x\| < \delta_0\}$, and using (3.1) we have that

$$\left[\sup_{\|t\| \leq C} \|S_n(t)\|\right]^2 \leq 2 a_n^{-4} C^2 \lambda_{\max}(M_{1n})$$

$$+ 2^{-1} a_n^{-6} p C^4 \sum_{i,j=1}^n \sum_{a=1}^p \lambda_{\max}(M_{2ni(a)}) \lambda_{\max}(M_{2nj(a)}).$$

Since $M_{1n}$ is nonnegative definite, its maximum eigenvalue can be dominated by its trace; from (3.19)–(3.21) it follows that $E[\sup_{\|t\| \leq C} \|S_n(t)\|]^2 = o(1)$. Note that

$$\inf_{|t|=C} \left\{ a_n^{-1} p^{-1/2} t^T \sum_{i=1}^n \phi_{ni}(\beta_0 + a_n^{-1} p^{1/2} t) \right\}$$

$$\geq -C \sup_{|t|=C} \|S_n(t)\| + C^2 l_{1n} - C a_n^{-1} p^{-1/2} \left\|\sum_{i=1}^n \phi_{ni}\right\|,$$

where $l_{1n} = \lambda_{\min}(2^{-1}(G_{1n} + G_{1n}^T))$, which from (3.22) is positive. Then by choosing $C$ large enough, we have that

$$\text{Prob}\left[\inf_{|t|=C}\left\{a_n^{-1}p^{-1/2}t^T\sum_{i=1}^n \phi_{ni}(\beta_0 + a_n^{-1}p^{1/2}t)\right\} > 0\right]$$

$$\geq 1 - \text{Prob}\left[a_n^{-1}p^{-1/2}\left\|\sum_{i=1}^n \phi_{ni}\right\| > Ck_2/2\right]$$

$$- \text{Prob}\left[\sup_{|t|=C}\|S_n(t)\| > Ck_2/2\right]$$

$$\geq 1 - \varepsilon \quad \text{for all } n \text{ sufficiently large.}$$

On the set where $\inf_{|t|=C} t^T \sum_{i=1}^n \phi_{ni}(\beta_0 + a_n^{-1}p^{1/2}t) > 0$, it then follows that $\sum_{i=1}^n \phi_{ni}(\beta_0 + a_n^{-1}p^{1/2}t) = 0$ for some $t \in \{t: \|t\| < C\}$ from continuity of $\phi_{ni}$'s and using Theorem 6.3.4 of Ortega and Rheinboldt (1970). Now (3.24) and (3.25) follow with a little algebra. The asymptotic normality is proved as in Theorem 3.1. □

Let $F_n(x) = \text{Prob}[s_n^{-1}c^T(\hat{\beta}_n - \beta_0) \leq x]$ and let $\Phi(\cdot)$ be the standard normal distribution function. Our model conditions are sufficient to argue that $\det(\phi_{1ni}(\hat{\beta}_n)) = 0$ has asymptotically negligible probability. In practice, this case is extremely unlikely. Hence define

$$F_{Bn}(x) = \mathbf{P}_B[\hat{s}_n^{-1}\sigma_n^{-1}c^T(\hat{\beta}_B - \hat{\beta}_n) \leq x]I_{\{\det(\phi_{1ni}(\hat{\beta}_n))\neq 0\}} + \Phi(x)I_{\{\det(\phi_{1ni}(\hat{\beta}_n))=0\}}$$

as the bootstrap distribution function estimator.

The next theorem is an analog of Theorem 3.2.

THEOREM 3.5. *Assume the conditions* (3.18)–(3.22) *and the bootstrap weights satisfy BW. There exists a sequence $\{\hat{\beta}_B\}$ of solutions of* (1.2) *such that if $p/a_n^2 \to 0$,*

$$\sigma_n^{-1}p^{-1/2}\sum_{i=1}^n \phi_{1ni}(\hat{\beta}_n)(\hat{\beta}_B - \hat{\beta}_n)$$

(3.26)

$$= -a_n^{-1}p^{-1/2}\sum_{i=1}^n W_i\phi_{ni}(\hat{\beta}_n) + r_{nB1},$$

*where $\|r_{nB1}\| = o_P(1)$. In addition, if* (3.23) *holds and BW and CLTW hold, then*

(3.27) $$\sup_x |F_{Bn}(x) - F_n(x)| \to 0 \quad \text{in probability.}$$

A sketch of the proof of this theorem is given after Remark 4.

REMARK 4. Lahiri (1992) has shown the consistency (and second-order accuracy) of an appropriate residual bootstrap for the usual $M$ estimation model with i.i.d. errors, known design and fixed $p$. Theorem 3.5 implies only the first-order consistency of the GBS and hence in particular of the PB, but for a much larger class of models.

In general GBS is not second-order accurate. First, $\hat{\beta}_n$ and $\hat{\beta}_B$ are biased for $\beta$ and $\hat{\beta}_n$, respectively, and the biases are not negligible in the second order. Further, as is known from the extensive literature on resampling, without an appropriate Studentization no resampling plan can hope to be second-order accurate. With appropriate bias correction and Studentization, the GBS can be made to be second-order accurate.

Define
$$\hat{g}_n^2 = n^{-1} \sum \phi_{ni}^2(\hat{\beta}_n) \quad \text{and} \quad \hat{g}_{nB}^2 = n^{-1} \sum W_i^2 \phi_{ni}^2(\hat{\beta}_n).$$

The following turns out to be the appropriate bias corrected Studentized statistic:
$$T_n = \hat{\gamma}_{1n} \hat{g}_n^{-1} [n^{1/2}(\hat{\beta}_n - \beta_0)] - 2^{-1} n^{-1/2} \hat{\gamma}_{1n}^{-2} \hat{g}_n^{-1} \hat{\gamma}_{2n} \mathbf{V}_{\text{GBS}},$$
$$T_{nB} = \hat{\gamma}_{1n} \hat{g}_{nB}^{-1} [\sigma_n^{-1} n^{1/2}(\hat{\beta}_B - \hat{\beta}_n)]$$
$$+ 2^{-1} n^{-1/2} \sigma_n \hat{g}_{nB}^{-1} \hat{\gamma}_{2n} [\sigma_n^{-1} n^{1/2}(\hat{\beta}_B - \hat{\beta}_n)]^2.$$

Chatterjee (1999) has shown that $T_{nB}$ is second-order accurate for $T_n$. However, there is ample scope for improvement on the conditions assumed there.

PROOF OF THEOREM 3.5. Let us concentrate on the set $\{\|\hat{\beta}_n - \beta_0\| < \delta_0/2\}$, since the complement of this set can be shown to have negligible contribution.

There we have
$$\mathbf{P}_B \left[ \sigma_n^{-1} p^{-1/2} a_n^{-1} \left\| \sum w_i \phi_{ni}(\hat{\beta}_n) \right\| > K \right]$$
$$\leq kK^{-2} p^{-1} a_n^{-2} \left[ \sum_{i=1}^n \|\phi_{ni}\|^2 + \|\hat{\beta}_n - \beta_0\|^2 \sum_{i=1}^n \sum_{a=1}^p \|\phi_{1nia}\|^2 \right.$$
(3.28)
$$\left. + \|\hat{\beta}_n - \beta_0\|^4 \sum_{i=1}^n \sum_{a=1}^p \lambda_{\max}^2(M_{2nia}) \right]$$
$$= K^{-2} O_P(1).$$

Thus for fixed $\delta_1, \delta_2 > 0$, by choosing $K$ large enough we have

(3.29) $\quad \text{Prob} \left[ \mathbf{P}_B \left[ \sigma_n^{-1} p^{-1/2} a_n^{-1} \left\| \sum w_i \phi_{ni}(\hat{\beta}_n) \right\| > K \right] > \delta_1 \right] < \delta_2.$

Let
$$S_{nB}(t) = \sigma_n^{-1} p^{-1/2} a_n^{-1} \sum_{i=1}^{n} w_i [\phi_{ni}(\hat{\beta}_n + \sigma_n p^{1/2} a_n^{-1} t) - \phi_{ni}(\hat{\beta}_n)] - \Gamma_{1n}^{T}(\hat{\beta}_n)t.$$

On the set $\{\|t\| \leq C\} \cap \{\|\hat{\beta}_n - \beta_0\| < \delta_0/2\}$, we have for large $n$

$$\|S_{nB}(t)\|^2 \leq 2\sigma_n^2 a_n^{-4} C^2 \lambda_{\max}\left(\sum_{a=1}^{p} \sum_{i,j=1}^{n} W_i W_j \phi_{1nia}(\hat{\beta}_n) \phi_{1nja}(\hat{\beta}_n)^T\right)$$
$$+ \sigma_n^2 p a_n^{-6} C^4 \sum_{a=1}^{p}\left(\sum_{i=1}^{n} w_i \lambda_{\max}(M_{2nia})\right)^2$$
$$= T_1 + T_2 \quad \text{say}.$$

With some algebra it can be shown that $\sum_{j=1}^{2} \mathbf{P}_{B}[T_j > K] = o_P(1)$, thus

(3.30)
$$\mathbf{P}_{B}\left[\sup_{\|t\| \leq C} \|S_{nB}(t)\| > 2K\right] \leq \sum_{j=1}^{2} \mathbf{P}_{B}[T_j > K] + O_P(a_n^{-1} p^{1/2})$$
$$= o_P(1).$$

Now observe that
$$\inf_{|t|=C}\left\{\sigma_n^{-1} a_n^{-1} p^{-1/2} t^T \sum_{i=1}^{n} w_i \phi_{ni}(\hat{\beta}_n + \sigma_n a_n^{-1} p^{1/2} t)\right\}$$
$$\geq -C \sup_{|t|=C} \|S_{nB}(t)\| + C^2 \hat{l}_{1n} - C\sigma_n^{-1} a_n^{-1} p^{-1/2} \left\|\sum_{i=1}^{n} w_i \phi_{ni}(\hat{\beta}_n)\right\|,$$

where $\hat{l}_{1n} = \lambda_{\min}(2^{-1}(\Gamma_{1n}(\hat{\beta}_n) + \Gamma_{1n}^{T}(\hat{\beta}_n)))$. Notice that $\hat{l}_{1n} > k_2/2$ with probability $1 - o(1)$, for the constant $k_2$ from (3.22). By choosing $C$ large enough, from (3.28), (3.29) and (3.30) we have that on the set $\{\|\hat{\beta}_n - \beta_0\| < \delta_0/2\}$,

$$\mathbf{P}_{B}\left[\inf_{|t|=C}\left\{\sigma_n^{-1} a_n^{-1} p^{-1/2} t^T \sum_{i=1}^{n} w_i \phi_{ni}(\hat{\beta}_n + \sigma_n a_n^{-1} p^{1/2} t)\right\} > 0\right]$$
$$\geq 1 - \mathbf{P}_{B}\left[\sigma_n^{-1} a_n^{-1} p^{-1/2} \left\|\sum_{i=1}^{n} w_i \phi_{ni}\right\| > Ck_2/2\right]$$
$$- \mathbf{P}_{B}\left[\sup_{|t|=C} \|S_{nB}(t)\| > Ck_2/2\right].$$

Thus for fixed $\delta_1, \delta_2 > 0$, we have that for $C$ large enough for all large $n$,

$$\text{Prob}\left[\mathbf{P}_{B}\left[\inf_{|t|=C}\left\{\sigma_n^{-1} a_n^{-1} p^{-1/2} t^T\right.\right.\right.$$

$$\times \sum_{i=1}^{n} w_i \phi_{ni}(\hat{\beta}_n + \sigma_n a_n^{-1} p^{1/2} t) \bigg\} > 0 \bigg] < 1 - \delta_1 \bigg]$$

$$\leq \text{Prob}\bigg[\mathbf{P}_{\mathrm{B}}\bigg[\sigma_n^{-1} a_n^{-1} p^{-1/2} \bigg\|\sum_{i=1}^{n} w_i \phi_{ni}\bigg\| > Ck_2/2\bigg] > \delta_1/2\bigg]$$

$$+ \text{Prob}\bigg[\mathbf{P}_{\mathrm{B}}\bigg[\sup_{|t|=C} \|S_n(t)\| > Ck_2/2\bigg] > \delta_1/2\bigg] + O(a_n^{-1} p^{1/2})$$

$$\leq \delta_2.$$

On the set $\inf_{|t|=C}\{\sigma_n^{-1} a_n^{-1} p^{-1/2} t^T \sum_{i=1}^{n} w_i \phi_{ni}(\hat{\beta}_n + \sigma_n a_n^{-1} p^{1/2} t)\} > 0$, using the continuity of $\sum_{i=1}^{n} w_i \phi_{ni}(\cdot)$ and Theorem 6.3.4 of Ortega and Rheinboldt (1970), we have that $\sum_{i=1}^{n} w_i \phi_{ni}(\hat{\beta}_n + \sigma_n a_n^{-1} p^{1/2} t) = 0$ has a root $T_n$ in $|t| \leq C$. Putting $\hat{\beta}_B = \hat{\beta}_n + \sigma_n a_n^{-1} p^{1/2} T_n$, we get a solution to (1.2) which satisfies, for fixed $\varepsilon, \delta > 0$, $\text{Prob}[\mathbf{P}_{\mathrm{B}}[\sigma_n^{-1} a_n p^{-1/2} \|\hat{\beta}_B - \hat{\beta}_n\| \leq C] < 1 - \varepsilon] < \delta$ for all $n$ large enough. Now notice that with this $C$ fixed, we have actually shown that with $t = T_n$

$$\sigma_n^{-1} a_n p^{-1/2} \Gamma_{1n}(\hat{\beta}_n)(\hat{\beta}_B - \hat{\beta}_n) = -a_n^{-1} p^{-1/2} \sum_{i=1}^{n} W_i \phi_{ni}(\hat{\beta}_n) + r_{nB1},$$

where $\|r_{nB1}\| = o_P(1)$. This shows (3.26). $\square$

**Acknowledgments.** We thank the past Editors Professor Hans Künsch and Professor Jon Wellner, an Associate Editor and the referees for their very careful reading of earlier versions of the manuscript and for suggesting many changes that we were glad to make. We thank Professors William Sudderth and Galin Jones for help with our English. The second author would like to thank Purdue University for hospitality during the preparation of an earlier version of this manuscript.

School of Statistics  
University of Minnesota  
313 Ford Hall, 224 Church Street S.E.  
Minneapolis, Minnesota 55455  
USA  
e-mail: chatterjee@stat.umn.edu

Statistics and Mathematics Unit  
Indian Statistical Institute  
203 B. T. Road  
Kolkata 700108  
India